\documentclass[english,leqno]{amsart}

\usepackage{amsmath}

\usepackage{amssymb}

\usepackage{amsbsy}

\usepackage{enumerate}

\usepackage[alphabetic]{amsrefs}
\usepackage{amsrefs}

\usepackage{hyperref}

\usepackage{mathrsfs} 

\usepackage{txfonts}  

\usepackage{frcursive}
\usepackage{calligra}
\usepackage[T1]{fontenc}

\usepackage{verbatim}

\usepackage{stmaryrd} 

\usepackage{fancyhdr}  

\usepackage{fullpage}

\usepackage{comment}

\usepackage{algorithm}

\usepackage{algpseudocode}

\usepackage{tikz}

\usetikzlibrary{arrows,positioning}

\usepackage{makecell}

\usetikzlibrary{datavisualization.formats.functions}

\usepackage{pgfplots}

\usepackage{graphicx,amssymb}

\usepackage{ifthen}

\usepackage{epsfig}

\newtheorem{theo}{Theorem}

\newtheorem{prop}{Proposition}

\newtheorem{lemm}{Lemma}

\newtheorem{coro}{Corollary}

\newtheorem{defi}{Definition}

\newtheorem{rema}{Remark}

\newtheorem{theoalpha}{Theorem}

\newtheorem*{thank}{Acknowledgments}

\renewcommand{\le}{\leqslant}
\renewcommand{\leq}{\leqslant}

\renewcommand{\ge}{\geqslant}
\renewcommand{\geq}{\geqslant}

\newcommand{\vep}{\varepsilon}

\newcommand{\ds}{\displaystyle}
\newcommand{\rest}{\restriction}

\newcommand{\abar}{\ensuremath{\bar{a}}}
\newcommand{\bbar}{\ensuremath{\bar{b}}}

\newcommand{\nbar}{\ensuremath{\bar{n}}}

\newcommand{\zbar}{\ensuremath{\bar{z}}}

\newcommand{\bC}{\ensuremath{\mathbb C}}
\newcommand{\bE}{\ensuremath{\mathbb E}}

\newcommand{\bN}{\ensuremath{\mathbb N}}

\newcommand{\bP}{\ensuremath{\mathbb P}}

\newcommand{\bR}{\ensuremath{\mathbb R}}

\newcommand{\cA}{\ensuremath{\mathcal A}}

\newcommand{\cdist}[1]{\mathsf{c}_{#1}}

\newcommand{\co}{\mathrm{c}_0}

\newcommand{\abs}[1]{\lvert #1\rvert}
\newcommand{\bigabs}[1]{\big\lvert #1\big\rvert}

\newcommand{\norm}[1]{\| #1\|}
\newcommand{\bnorm}[1]{\Big\| #1\Big\|}

\newcommand{\esp}{\mathbb{E}}

\newcommand{\sgn}{\operatorname{sgn}}

\newcommand{\dtreek}{(x_{t})_{t \in D^{\leq k}\setminus\{\emptyset\}}}

\begin{document}

\allowdisplaybreaks

\title{On asymptotic B-convexity and infratype}

\author{F.~Baudier}
\address{F.~Baudier, Department of Mathematics, Texas A\&M University, College Station, TX 77843, USA}
\email{florent@tamu.edu}

\author{A.~ Fovelle}
\address{Institute of Mathematics (IMAG) and Department of Mathematical Analysis, University of Granada, 18071, Granada, Spain}
\email{audrey.fovelle@ugr.es}

\thanks{F. Baudier was partially supported by the National Science Foundation under Grant Number DMS-1800322 and DMS-20556040. A. Fovelle was partially supported by MCIN/AEI/10.13039/501100011033 grant PID2021-122126NB-C31, by ``Maria de Maeztu'' Excellence Unit IMAG, reference CEX2020-001105-M funded by MCIN/AEI/10.13039/501100011033, and by NSF Grant DMS-1800322.}
	


\begin{abstract}
In this note, we introduce and study the notions of asymptotic B-convexity and asymptotic infratype $p$, and we prove asymptotic analogs of a series of results due to Giesy \cite{Giesy66} and Pisier \cite{Pisier74}. In particular, we give a simplified proof of an asymptotic version of Pisier's $\ell_1$-theorem that was originally proven by Causey, Draga, and Kochanek in \cite{CDK19}. 
We also briefly discuss the notion of asymptotic stable type.
\end{abstract}

\maketitle

\setcounter{tocdepth}{3}

\section{Introduction}

First introduced by Beck in 1962 \cite{Beck62} to characterize the Banach spaces for which a certain strong law of large numbers for Banach-valued random variables holds, B-convexity\footnote{This property was aptly named property $(B)$ in \cite{Beck62} (where it was compared to a property $(A)$) and it seems that Giesy coined the term $B$-convexity in \cite{Giesy66}. James also uses the terminology uniformly non-$\ell_1^k$ for a space that is B-$(k,\delta)$-convex for some $\delta\in(0,1)$, as well as uniformly non-square when $k=2$ and uniformly non-octahedral when $k=3$.} was defined as follows: a Banach space $X$ is B-convex if there exist $k \in \bN$ and $\delta\in(0,1)$ such that for every finite sequence of vectors $x_1, \dots, x_k$ in the unit ball of $X$, denoted by $B_X$, there is a combination of signs $(\vep_1,\dots,\vep_k)\in \{-1,1\}^k$ such that
\begin{equation}
\label{eq:B-convex}
    \bnorm{\sum_{i=1}^k \vep_i x_i} \leq k(1- \delta).
\end{equation}
If \eqref{eq:B-convex} holds for $k \in \bN$ and $\delta\in(0,1)$, we say that $X$ is B-$(k, \delta)$-convex.

Soon after, Giesy \cite{Giesy66} carried out a thorough study of $B$-convexity and proved many of its basic properties. In particular, Giesy showed that $X$ is B-convex if and only if $\ell_1$ is not finitely representable in $X$ (or equivalently, $X$ does not contain the $\ell_1^k$'s uniformly), i.e., one cannot find a constant $\lambda\in[1,\infty)$ such that, for all $k \in \bN$, one can find $x_1, \dots, x_k \in S_X$, the unit sphere of $X$, satisfying $\lambda^{-1} \norm{\abar}_{\ell_1^k} \leq \norm{\sum_{i=1}^k a_i x_i } \leq  \norm{\abar}_{\ell_1^k}$ for every $\abar:=(a_1, \dots, a_k) \in \bR^k$. Giesy also showed that $B$-convexity is stable under linear isomorphism, and under passing to quotients or taking duals. While a B-convex Banach space with an unconditional basis is reflexive (see \cite{Giesy66} or \cite{James64}), James \cite{James74} constructed an example of a B-convex Banach space (in fact uniformly non-octahedral) that is not reflexive.

Finite representability of $\ell_1$, and thus $B$-convexity, is closely connected to the notion of Rademacher type. Recall that a Banach space has (Rademacher) type $p\in[1,2]$ if there is a constant $T_R\in(0,\infty)$ such that for every $k\in\bN$ and $x_1,\dots, x_k$ in $X$, we have
\begin{equation}
\label{eq:type}
\bE_{\vep\in \{\pm 1\}^k} \bnorm{ \sum_{i=1}^k \vep_i x_i }^p \leq T_R^p \sum_{i=1}^k \norm{x_i}^p.
\end{equation}
Note that every Banach space has type $1$. We say that a Banach space has non-trivial type if it has type $p$ for some $p>1$. It is easy to see that if $X$ has type $p$ with $p>1$, then $X$ does not contain the $\ell_1^k$'s uniformly. Remarkably, Pisier \cite{Pisier74} showed that the converse holds.
Pisier's proof goes via the notion of infratype $p$. A Banach $X$ has \emph{infratype $p \in [1,2]$} if there is a constant $T_I\in(0,\infty)$ such that for every $k\in\bN$ and $x_1,\dots, x_k$ in $X$, we have
\begin{equation}
\label{eq:infratype}
\min_{\vep\in \{\pm 1\}^k} \bnorm{ \sum_{i=1}^k \vep_i x_i }^p \leq T_I^p \sum_{i=1}^k \norm{x_i}^p.
\end{equation}
It is plain that type $p$ implies infratype $p$. Several sequences of parameters, measuring in some sense the B-convexity, type $p$, or infratype $p$ behavior of a Banach space were introduced in \cite{Pisier74}. The submultiplicativity of certain of these parameter sequences was crucial in showing that a Banach space $X$ is B-convex if and only if $X$ has non-trivial infratype if and only if $X$ has non-trivial type if and only if $X$ does not contain the $\ell_1^k$'s uniformly.

Beauzamy \cite{Beauzamy76-type} proved an operator version of Pisier's $\ell_1$-theorem. He showed that an operator $T\colon X\to Y$ has (Rademacher) subtype if and only if $\ell_1$ is not finitely represented in $T$ (we refer to \cite{Beauzamy76-type} for the precise definitions that are not needed in this note). The spatial result of Pisier can be recovered by considering the identity operator and showing that the identity operator on a Banach space $X$ has subtype if and only if $X$ has non-trivial type (but this latter fact is not necessarily true for arbitrary operators). More recently, Causey, Draga, and Kochanek \cite{CDK19} proved an asymptotic version of Beauzamy's result. The asymptotic approach of \cite{CDK19} is from the point of view of asymptotic structure as laid out by Maurey, Milman, and Tomczak-Jaegermann in \cite{MMTJ95}. In \cite{CDK19}, the authors develop a concept of block structures for which they define various notions of asymptotic type. Then, they extend the spatial notion of asymptotic structure of \cite{MMTJ95} to the operator setting by defining the notion of asymptotic structure of an operator using their framework of block structures, via certain block structures that are induced by operators. This allows them to define various asymptotic notions of type for operators for which they prove asymptotic analogs of Beauzamy's results. By specializing to the identity operator on a Banach space $X$ they obtain various asymptotic notions of type for $X$ and derive that a Banach space has non-trivial asymptotic type if and only if $\ell_1$ is not asymptotically finitely representable in $X$. This last concept is also introduced in \cite{CDK19} via the block structures approach. These notions can be reformulated in terms of properties of branches of certain weakly null trees of finite but arbitrarily large height as shown in \cite{CDK19}. As the authors in \cite{CDK19} are mainly interested in operator ideals, they mostly follow Beauzamy's approach adjusted to the framework of block structures. Beauzamy's techniques are inspired by Pisier's original work but the operator setting being slightly different the arguments are also slightly different.  
In this note, we take an approach opposite to \cite{CDK19}. Firstly, we focus on spatial results. Secondly, our starting point is the definition of asymptotic notions in terms of properties of branches of weakly null trees of finite but arbitrarily large height thus avoiding the (rather tedious) notion of block structures altogether. Lastly, we follow Pisier's original approach. This requires introducing notions that were not considered in \cite{CDK19}, and we prove asymptotic analogs of Pisier's results from \cite{Pisier74} and the needed results from \cite{Giesy66} about $B$-convexity. 

In Section \ref{sec:B-convexity}, we introduce the notion of asymptotic B-convexity and prove the following equivalences due to Giesy in the local setting. 

\begin{theoalpha} 
\label{thmA:asymp-B-convex-ell1}
Let $X$ be a Banach space. The following assertions are equivalent.
\begin{enumerate}[(i)]
\item $X$ is not asymptotically B-convex.
\item $\ell_1$ is asymptotically finitely representable in $X$.
\item $\ell_1$ is crudely asymptotically finitely representable in $X$.
\end{enumerate}
\end{theoalpha}

In Section \ref{sec:infratype}, we introduce the notion of asymptotic infratype and prove the following theorem.  

\begin{theoalpha} 
\label{thmB:asymp-infratype<->asymp-type}
Let $X$ be a Banach space. The following assertions are equivalent.
\begin{enumerate}[(i)]
\item $X$ has non-trivial asymptotic infratype.
\item $X$ has non-trivial asymptotic type.
\item $X$ is asymptotically B-convex.
\end{enumerate}
\end{theoalpha}

Almost all the proofs of the local results discussed above rely on introducing certain sequences of parameters and leveraging the fact that some of these sequences are submultiplicative. The submultiplicative proofs in the local setting are completely elementary. However, when considering the corresponding sequences of parameters for the asymptotic notions defined in terms of weakly null trees, these proofs become significantly more delicate and form the main contribution of this note. By way of comparison, the submultiplicative proofs in the block structures approach in \cite{CDK19} remain completely elementary, and the technical burden is encompassed in the definition of the block structures as well as showing that the asymptotic notions defined in terms of block structures are equivalent to the asymptotic notions defined in terms of trees.

Moreover, Pisier's original approach in the spatial setting used three sequences of parameters while Beauzamy only used two sequences in the operator setting. For our asymptotic results in the spatial setting, we only use two sequences of parameters as well and simplify further the argument by using the contraction principle (that is neither used in \cite{Pisier74} nor \cite{Beauzamy76-type} for this purpose).
We believe that our presentation gives the most streamlined approach to prove the equivalence above in either the local or asymptotic setting, and to make the note essentially self-contained we provide all the details even for arguments similar to their local counterparts.

Finally, for completeness, we briefly discuss the notion of asymptotic stable type in Section \ref{sec:stable-type}.

\section{Asymptotic type and asymptotic finite representability}
\label{sec:type-fr}
In this section, we recall the tree-versions of the notions of asymptotic type and asymptotic finite representability introduced in \cite{CDK19}. Since we want to study these asymptotic notions in arbitrary Banach spaces (and not only those with separable dual), we first recall some notions related to trees indexed by directed sets. Let $D$ be an arbitrary set. We denote by $\emptyset$, the empty sequence, and $D^{\le k}=\{\emptyset\}\cup \cup_{i=1}^k D^i$, the set of all finite sequences in $D$ of length at most $n$.
For $s,t\in D^{\le k}$, we let $s\frown t$ denote the concatenation of $s$ with $t$, and we simply write $s\frown a$ instead of $s\frown (a)$, whenever $a\in D$. 
We let $\abs{t}$ denote the length of $t\in D^{\le k}$, i.e., $\abs{t}=i$ if and only if $t\in D^i$. For $0\le i\le\abs{t}$, we let $t_{\rest i}$ denote the initial segment of $t$ having length $i$, where $t_{\rest 0}=\emptyset$. This means that if $t=(t_1,\dots, t_j)\in D^{\le k}$ and $0\le i\le j$, then $t_{\rest i} = (t_1,\dots, t_i)$. If $s\in D^{\le k}$, we let $s\prec t$ denote the relation that $s$ is a proper initial segment of $t$, i.e., $s=t_{\rest i}$ for some $0\le i<\abs{t}$. We denote $s \preceq t$ if $s\prec t$ or $s=t$.  It is well known that  $(D^{\le k}, \prec)$ is a tree in the set-theoretic sense.

If $(D, \leq)$ is a directed set and $(x_t)_{t\in D^{\le k}}\subset X$, we say that $(x_t)_{t\in D^{\le k}}$ is  a \emph{weakly null tree} in $X$, provided that for each $t\in D^{\le k -1}$, $(x_{t\frown s})_{s\in D}$ is a weakly null net. A function $\varphi \colon D^{\le k}\to D^{\le k}$ is said to be a \emph{pruning} provided that 
\begin{enumerate}[(i)]
\item $\varphi$ preserves the tree ordering, i.e., if $s\prec t$, then $\varphi(s)\prec \varphi(t)$,
\item $\varphi$ preserves the length, i.e., $|\varphi(t)|=|t|$ for all $t\in D^{\le k}$, and moreover 
\item if $\varphi((t_1, \ldots, t_j))=(s_1, \ldots, s_j)$, then $t_i\leqslant s_i$ for all $1\le i \le j\le k$. 
\end{enumerate} 
One reason one might want to prune a tree is to stabilize maps that are defined on the leaves of trees. The first lemma is a Ramsey-type theorem in the context of trees on directed sets.

\begin{lemm}
\label{lem:pruning1}
Let $(D,\leq)$ be a directed set, $F$ a finite set, $k\in \bN$, and $f \colon D^k \to F$ a function. Then, there exists a pruning $\varphi \colon D^{\le k}\to D^{\le k}$ such that $f \circ \varphi$ is constant on $D^k$.
\end{lemm}


The definition below, of a notion of asymptotic type from the point of view of asymptotic structure, can be extracted from \cite[Definition 5.7 and Lemma 6.2]{CDK19}.

\begin{defi}
\label{def:asymptotic-type}
    A Banach space $X$ has \emph{asymptotic (Rademacher) type $p$} if and only if there exists a constant $T>0$ such that for every directed set $D$, every weakly null tree $\dtreek \subset B_X$, and every scalar sequence $(a_1,\dots, a_k)\in \bR^k$, there is $t\in D^k$ such that  
\begin{equation}
\label{eq:asymptotic-type-p}
\bE_\vep \bnorm{\sum_{j=1}^k \vep_j a_j x_{t_{\rest j} } }^p \le T^p \sum_{j=1}^k \abs{a_j}^p.
\end{equation}
\end{defi}

It is easy to see that asymptotic type $p$ is implied by type $p$ and $p$-asymptotic uniform smoothness.
The existence of a branch $t$ satisfying \eqref{eq:asymptotic-type-p} implies the existence of a pruning of $D^{\le k}$ so that \eqref{eq:asymptotic-type-p} holds for every branch in the pruned tree.

\begin{prop}
$X$ has asymptotic Rademacher type $p$ if and only if there exists $T>0$ such that, for every directed set $D$, every $k \in \bN$ and every weakly null tree $(x_t)_{t \in D^{\leq k}\setminus\{\emptyset\}} \subset B_X$, there exists a pruning $\varphi \colon D^{\leq k}\to D^{\le k}$ so that for every scalar sequence $(a_1,\dots, a_k)\in \bR^k$ and $t \in D^k$, 
\begin{equation}
    \bE_\vep\bnorm{\sum_{j=1}^k \vep_j a_j x_{\varphi(t)_{\rest j}} }^p \leq T^p \sum_{j=1}^k \abs{a_j}^p.
\end{equation}
\end{prop}

\begin{proof}
For the non-trivial implication, assume that $X$ has asymptotic Rademacher type $p$ with constant $T>0$. It follows from Lemma \ref{lem:pruning1}, that for every directed set $D$, every $k \in \bN$, every weakly null tree $(x_t)_{t \in D^{\leq k}\setminus\{\emptyset\}} \subset B_X$, and every scalar sequence $\abar:=(a_1,\dots, a_k)\in S_{\ell_p^k}$, there exists a pruning $\varphi \colon D^{\leq k}\to D^{\leq k}$ so that for all $t\in D^k$, 
\begin{equation}
     \bE_\vep \bnorm{\sum_{j=1}^k \vep_j a_j x_{\varphi(t)_{\rest j}} }^p \leq T^p. 
\end{equation}
Let $N$ be a $\delta$-net of $S_{\ell_p^k}$. Since $N$ is finite, we can find a pruning $\phi \colon D^{\leq k}\to D^{\leq k}$ so that for all $\bbar \in N$ and $t \in D^k$, 
\begin{equation}
     \bE_\vep \bnorm{\sum_{j=1}^k \vep_j b_j x_{(\varphi\circ\phi)(t)_{\rest j}} }^p \leq T^p. 
\end{equation}
Fix $t \in D^k$ and $\abar \in S_{\ell_p^k}$, and pick $\bbar \in N$ so that $\bnorm{\abar -\bbar}_{\ell_p^k} \leq \delta$. Then
\begin{align*}
\bE_\vep \bnorm{\sum_{j=1}^k \vep_j a_j x_{(\varphi\circ\phi)(t)_{\rest j}} }^p & \le 2^{p-1} (\bE_\vep \bnorm{ \sum_{j=1}^k \vep_j b_j x_{(\varphi\circ\phi)(t)_{\rest j}} }^p + \bE_\vep \bnorm{\sum_{j=1}^k \vep_j (a_j-b_j) x_{(\varphi\circ\phi)(t)_{\rest j}} }^p) \\
& \leq 2^{p-1}T^p + 2^{p-1}(k\delta)^{p}\\
& \leq 2^pT^p,
\end{align*}
if $\delta$ is chosen small enough. We conclude by homogeneity.
\end{proof}

The following definition was also introduced in \cite{CDK19}. 

\begin{defi}
    Let $p\in[1,\infty]$ and $C\ge 1$. We say that $\ell_p$ is $C$-crudely asymptotically finitely representable in $X$ if there are constants $A,B>0$ with $AB\le C$ such that for all $k\in \bN$, there are a directed set $D$ and a weakly null tree $(x_t)_{t \in D^{\leq k}\setminus\{\emptyset\}} \subset S_X$, such that for every scalar sequence $(a_1,\dots, a_k)\in \bR^k$ and every $t\in D^k$,
    \begin{equation*}
    \frac{1}{A}(\sum_{j=1}^k \abs{a_j}^p)^{1/p} \le  \bnorm{ \sum_{j=1}^k a_j x_{t_{\rest j}} }  \le B(\sum_{j=1}^k \abs{a_j}^p)^{1/p},
\end{equation*}
with the obvious convention when $p=\infty$.
    We say that $\ell_p$ is \emph{asymptotically finitely representable in $X$} if it is $C$-crudely asymptotically finitely representable in $X$ for every $C>1$.
\end{defi}

\section{Asymptotic B-convexity}
\label{sec:B-convexity}

The following definition of an asymptotic notion of B-convexity is natural from the point of view of asymptotic structure.

\begin{defi}
\label{defi:asymp-B-convexity}
Let $k \in \bN$ and $\delta\in(0,1)$. We say that a Banach space $X$ is \emph{asymptotically B-$(k, \delta)$-convex} if for every directed set $D$, every weakly null tree $\dtreek \subset B_X$, there exist $t \in D^k$ and an arrangement of signs $(\vep_j)_{j=1}^k \in \{-1,1\}^k$ such that
\begin{equation}
\bnorm{\sum_{j=1}^k \vep_j x_{t_{\rest j}}} \leq k-\delta.
\end{equation}
We say that $X$ is \emph{asymptotically B-convex} if it is asymptotically B-$(k, \delta)$-convex for some $k \in \bN$ and $\delta > 0$.
\end{defi}

The next proposition is the asymptotic analog of the fact that if $X$ is not asymptotically B-$(k,\delta)$ convex, then $\ell_1^k$ is $\frac{1}{1-\delta}$-crudely finitely representable in $X$.

\begin{prop}
\label{pro:not-asymp-B-convex->ell1} 
Let $k \in \bN$ and $\delta\in(0,1)$.
If $X$ is not asymptotically $(k,\delta)$-convex, then $\ell_1^k$ is $\frac{1}{1-\delta}$-crudely asymptotically finitely representable in $X$.
\end{prop}

\begin{proof}
If $X$ is not asymptotically B-$(k,\delta)$ convex, then there exist a directed set $D$ and a weakly null tree $\dtreek \subset B_X$ such that for all $t \in D^k$, 
\begin{equation}
\label{eq:not-asymp-B-convex->ell1}
 \inf_{(\vep_j)_{j=1}^k \in \{\pm 1 \}^k} \bnorm{\sum_{j=1}^k \vep_j x_{t_{\rest j}} } \geq k-\delta.
 \end{equation}
Let $(a_j)_{j=1}^k\in \bR^k$ such that $\sum_{j=1}^k \abs{a_j}=1$ and $t \in D^k$. If for every $1 \leq j \leq k$, we let $\vep_j:=\sgn(a_j)$ so that $a_j=\vep_j\abs{a_j}$, then 
\begin{equation*}
\bnorm{ \sum_{j=1}^k a_j x_{t_{\rest j}} } \ge \bnorm{ \sum_{j=1}^k \vep_j x_{t_{\rest j}}} -\bnorm{ \sum_{j=1}^k \vep_j(\abs{a_j}-1) x_{t_{\rest j}} } \stackrel{\eqref{eq:not-asymp-B-convex->ell1}}{\ge} k-\delta -\sum_{j=1}^k (1- \abs{a_j}) \ge 1-\delta,    
\end{equation*}
and hence by homogeneity we have that for all $(a_j)_{j=1}^k\in \bR^k$ and all $t\in D^k$,
\begin{equation*}
    (1-\delta)\sum_{j=1}^k \abs{a_j} \le  \bnorm{ \sum_{j=1}^k a_j x_{t_{\rest j}} }  \le \sum_{j=1}^k \abs{a_j}.
\end{equation*}
\end{proof}

When $\ell_1$ is asymptotically finitely representable in $X$ it is easy to see that $X$ cannot be asymptotically B-convex, and hence the following result, due to Giesy \cite{Giesy66} in the local setting, follows immediately from Propostion \ref{pro:not-asymp-B-convex->ell1}.
\begin{prop}
\label{pro:asymp-B-convex-ell1}
 $X$ is not asymptotically B-convex if and only if  $\ell_1$ is asymptotically finitely representable in $X$.
\end{prop}

It is possible to relax the condition of asymptotic finite representability to crude asymptotic finite representability. This can be done by proving an asymptotic analog of James $\ell_1$-distortion theorem, but it will also be a consequence of the submultiplicative nature of a natural parameter that can be associated with any Banach space and detects its asymptotic B-convexity. This approach was used by Pisier \cite{Pisier74} to prove the local analog of this result. We develop this approach in the asymptotic setting as it will be useful for other purposes in the sequel.

For every $k \in \bN$, let us denote by $\bbar_X(k)$ the smallest possible constant $b \geq 0$ such that for every directed set $D$ and every weakly null tree $\dtreek \subset B_X$, there exists $t \in D^k$ such that
\begin{equation}
 \min_{(\vep_j)_{j=1}^k \in \{\pm 1\}^k} \bnorm{\sum_{j=1}^k \vep_j x_{t_{\rest j}} } \leq b k .
 \end{equation}

Note that $\bbar_X(1)=1$. If $\ell_1$ is crudely asymptotically finitely representable in $X$, it is easy to see that $\inf_{k\ge 1} \bbar_X(k)>0$. In order to prove that if $X$ is asymptotically B-convex, then $\ell_1$ is not crudely asymptotically finitely representable in $X$, it is then enough to prove that $\liminf\limits_{k \to \infty} \bbar_X(k) = 0$ whenever $X$ is asymptotically B-convex. This fact will be a consequence of the submultiplicativity of the sequence $(\bbar_X(k))_{k\in \bN}$. The proof of this submultiplicativity property in the asymptotic setting is arguably more technical than in the local setting. 

In the next lemma, we also record useful properties of the sequence $(\bbar_X(k))_{k\in \bN}$.

\begin{lemm} Let $X$ be an infinite-dimensional Banach space. Then,
\label{lem:b-submul}
\begin{enumerate}[(i)]
\item $X$ fails to be asymptotically B-convex if and only if $\bbar_X(k)=1$ for all $k\ge 1$.
\item The map $k\mapsto k\bbar_X(k)$ is non-decreasing.
\item For all $k\in \bN$, $\frac{1}{k} \le \bbar_X(k) \le 1$.
\item For all $k,l \in \bN$, $\bbar_{X}(kl)\le \bbar_X(k)\bbar_X(l)$.
\end{enumerate}
\end{lemm}

\begin{proof}
Assertion $(i)$ is trivial and assertion $(iii)$ follows from assertion $(ii)$.
The proof of assertion $(ii)$ is standard and elementary. Indeed, given $k\in \bN$, $D$ be a directed set, and $\dtreek$ be a weakly null tree in $B_X$. Let $(y_{t})_{t \in D^{\leq k+1}\setminus\{\emptyset\}}$ defined by $y(t)=x(t)$ if $\abs{t}\le k$ and $y_t=0$ otherwise. Then $(y_{t})_{t \in D^{\leq k+1}\setminus\{\emptyset\}}$ is a weakly null tree in $B_X$ and there is $t\in D^{\le k+1}$ such that
\begin{equation*}
 \min_{(\vep_j)_{j=1}^k \in \{\pm 1\}^k} \bnorm{\sum_{j=1}^k \vep_j x_{t_{\rest j}} } =  \min_{(\vep_j)_{j=1}^{k+1} \in \{\pm 1\}^k} \bnorm{\sum_{j=1}^{k+1} \vep_j y_{t_{\rest j}} } \leq \bbar_X(k+1) (k+1) = \bbar_X(k+1)\frac{k+1}{k}k, 
 \end{equation*}
and hence $\bbar_X(k)\le\bbar_X(k+1)\frac{k+1}{k}$, and the claim is proved.
Assertion $(iv)$ is the main assertion. Let $k,l \in \bN$, $D$ be a directed set, and $(x_{t})_{t \in D^{\leq kl}\setminus\{\emptyset\}} \subset B_X$ be a weakly null tree. The first order of business is to prune the tree so that certain properties are true for all branches in the pruned tree. It will help to think of a tree of height $kl$ as a tree of height $l$ where on top of each of its leaves sits another tree of height $l$, and so on for $k$ steps. Pick $b_k > \bbar_X(k)$, $b_l > \bbar_X(l)$ and consider the set 
\begin{equation*}
    \cA:=\Big\{ t \in D^{kl} \colon \max_{0\le i< k} \min_{\vep\in \{\pm 1\}^l} \bnorm{\sum_{j=1}^l \vep_j x_{t_{\rest il+j}} } \le lb_l
 \Big\}.
\end{equation*}
 Since a pruned weakly null tree is still weakly null and $b_l > \bbar_X(l)$, for any pruning $\varphi \colon D^{\le kl} \to D^{\le kl}$ one can find $t\in D^{kl}$ such that $\varphi(t)\in \cA$. 
Then, it follows from Lemma \ref{lem:pruning1}, by considering the indicator function of $\cA$, that there is a pruning $\phi\colon D^{\le kl} \to D^{\le kl}$ such that for all $t \in D^{kl}$ and all $0\le i<k$,
\begin{equation*}
\min_{\vep\in \{\pm 1\}^l} \bnorm{\sum_{j=1}^l \vep_j x_{\phi(t)_{\rest il+j}} }\le lb_l.
\end{equation*}
Invoking Lemma \ref{lem:pruning1}, for $F=\{-1,1\}^{kl}$ which is a finite set, there is a combination of signs $\bar{\vep}\in \{ \pm 1 \}^{kl}$ and a further pruning, still denoted $\phi$, such that for all $t \in D^{kl}$ and all $0\le i<k$, 
\begin{equation}
\label{eq:eq1}
\bnorm{\sum_{j=1}^l \vep_{li+j} x_{\phi(t)_{\rest il+j}} }\le lb_l.
\end{equation}

Now, we construct a weakly null tree $(y_{u})_{u\in D_X^{\le k}\setminus\{\emptyset\}}$, for some directed set $D_X$, and for all $u\in D_X^{\le k}\setminus\{\emptyset\}$, $t^{(u_1)}\prec t^{(u_1,\dots,u_i)} \prec t^{(u_1,\dots,u_{k})} \in D^{\le kl}\setminus\{\emptyset\}$ 
 with $t^{(u_1,\dots,u_i)}\in D^{il}$, such that $y_{(u_1,\dots,u_i)}= \frac{1}{l b_l} \sum_{j=1}^l \vep_{(i-1)l+j} x_{\phi(t^{(u_1,u_2, \dots,u_i)})_{\rest (i-1)l+j}}$. We take $D_X$ to be a weak convex neighborhood basis of $0$ in $X$, directed by reverse inclusion. 
Let $U \in D_X$. Since the pruned tree is weakly null we can certainly find $t^{(U)}\in D^l$ such that $\frac{1}{b_l} \vep_j x_{\phi(t^{(U)})_{\rest j}} \in U$ for all $1 \leq j \leq l$. Then we let $y_U := \frac{1}{l b_l} \sum_{j=1}^l \vep_j x_{\phi(t^{(U)})_{\rest j}}$. It follows from the convexity of the neighborhoods that $y_U \in U$.
The construction of $(y_U)_{U \in D_X}$ clearly makes it weakly null and because of $\eqref{eq:eq1}$ it is in $B_X$.

Now, fix $(U_1,U) \in D_X^2$. Since the tree is weakly null we can certainly find $t^{(U_1,U)}\in D^{2l}$ with $t^{U_1}\prec t^{(U_1,U)}$ and such that $\frac{1}{b_l} \vep_{l+j} x_{\phi(t^{(U_1,U)})_{\rest l+j}} \in U$ for all $1 \leq j \leq l$. Then we let $y_{(U_1,U)} := \frac{1}{l b_l} \sum_{j=1}^l \vep_{l+j} x_{\phi(t^{(U_1,U)})_{\rest l+j}} \in  U$, and yet again $(y_{(U_1,U)})_{U \in D_X}$ is weakly null and in $B_X$ for all $U_1\in D_X$.

Iterating we can find for all $1\le i \le k$ and $(U_1,U_2,\dots, U_i)\in D_X^i$, a node $t^{(U_1,U_2, \dots,U_i)}\in D^{il}$ such that 
\begin{enumerate}
    \item $\frac{1}{b_l} \vep_{(i-1)l+j} x_{\phi(t^{(U_1,\dots,U_i)})_{\rest (i-1)l+j} } \in U_i$ for all $1 \leq j \leq l$,
    \item $t^{(U_1,\dots,U_{i-1})}\prec t^{(U_1,\dots,U_i)}$.
\end{enumerate}
If for all $1\le i \le k$ we define 
\begin{equation}
y_{(U_1,\dots,U_i)}:= \frac{1}{l b_l} \sum_{j=1}^l \vep_{(i-1)l+j} x_{ \phi(t^{(U_1,U_2, \dots,U_i)})_{\rest (i-1)l+j}},
\end{equation}
then the tree $(y_{u})_{u \in D_X^{\le k}\setminus \{\emptyset\}}$ is weakly null and in $B_X$. Now, since $b_k> \bbar_X(k)$, we can find $u:=(U_1,\dots,U_k) \in D_X^k$ and $\delta_1, \dots, \delta_k \in \{-1,1\}$ such that 
\begin{equation}
\bnorm{ \sum_{i=1}^k \delta_i y_{u_{\rest i}}} \leq k b_k .
\end{equation}
This means that
\begin{equation}
\bnorm{ \sum_{i=1}^k \delta_i \frac{1}{lb_l}\sum_{j=1}^l \vep_{(i-1)l+j} x_{\phi(t^{(U_1,U_2, \dots,U_i)})_{\rest (i-1)l+j}} } \leq k b_k,
\end{equation}
and hence 
\begin{equation}
\bnorm{ \sum_{i=1}^k \sum_{j=1}^l \delta_i\vep_{(i-1)l+j} x_{\phi(t^{(U_1,U_2, \dots,U_i)})_{\rest (i-1)l+j}} } \leq b_lb_k lk ,
\end{equation}
We have thus extracted a branch $\phi(t^{u})\in D^{kl}$ and a choice of signs $\bar{\eta}\in \{-1,1\}^{kl}$, namely 
$$\bar{\eta}:=(\delta_1\vep_1, \dots,\delta_1\vep_l, \delta_2\vep_{l+1}, \dots, \delta_2\vep_{2l},\dots, \delta_k\vep_{(k-1)l+1},\dots, \delta_k\vep_{kl}),$$ such that 
$\bnorm{\sum_{r=1}^{kl} \eta_r x_{\phi(t^{u})_{\rest r} } } \le b_lb_k lk$, from which it follows that $\bbar_X(kl)\le \bbar_X(l)\bbar_X(k)$.

\end{proof}

\begin{rema}
In the local setting, $b_X(k)$ the local analog of $\bbar_X(k)$ satisfies $b_X(k)\ge \frac{1}{\sqrt{k}}$; this is a consequence of Dvoretzky-Rogers theorem. This is no longer true in the asymptotic setting since, for instance, $\bbar_{\co}(k)=\frac1k$.
This last observation, also shows that while $\ell_1$ is certainly finitely representable in $\co$, it is not asymptotically finitely representable in $\co$.
\end{rema}

The following corollary, which contains Theorem A, follows from Proposition \ref{pro:asymp-B-convex-ell1} and the discussion above. In particular, it shows that asymptotic B-convexity is an isomorphic invariant.
\begin{coro} Let $X$ be a Banach space. The following assertions are equivalent.
\label{cor:asymp-B-convex-ell1}
\begin{enumerate}[(i)]
\item $X$ is not asymptotically B-convex.
\item $\ell_1$ is asymptotically finitely representable in $X$.
\item $\ell_1$ is crudely asymptotically finitely representable in $X$.
\item $\inf_{k\ge 1}\bbar_X(k)>0$.
\end{enumerate}
\end{coro}

\begin{proof}
Only $(iv) \implies (i)$ is missing, but if $X$ is asymptotically B-convex then there is $k_0\ge 2$ and $\delta\in(0,1)$ such that $X$ is asymptotically B-$(k_0,\delta)$-convex. This means that $\bbar_X(k_0)<1$ and by submultiplicativity $\lim\limits_{n\to \infty}\bbar_X(k_0^n)=0$.
\end{proof}

Since it was shown in \cite{CDK19} that the property ``$\ell_1$ is not asymptotically finitely representable in'' is a three-space property we have the following corollary.

\begin{coro}
 Asymptotic B-convexity is a three-space property, i.e., if $Y$ is a closed subspace of $X$ and both $Y$ and $X/Y$ are asymptotically B-convex, then $X$ is asymptotically B-convex.
\end{coro}

In the local setting, it is well known that if a Banach space admits an equivalent uniformly smooth norm, then $X$ admits an equivalent norm that is uniformly non-square, i.e., B-$(2,\delta)$-convex for some $\delta>0$. A similar result holds for the corresponding asymptotic notions.

\begin{prop} 
\label{prop:AUS->B-convex}
If $X$ admits an equivalent asymptotic uniformly smooth norm, then $X$ admits an equivalent norm that is asymptotically B-$(2,\delta)$-convex for some $\delta>0$.
\end{prop}

\begin{proof}
If $X$ admits an equivalent AUS norm, then it admits an equivalent norm and $p\in (1, \infty)$ such that
\begin{equation}
\limsup_\lambda \norm{ x + x_\lambda }^p \leq \norm{ x }^p + \limsup_\lambda \norm{ x_\lambda }^p,
\end{equation}
whenever $x \in X$ and $(x_\lambda)$ is a bounded weakly null net in $X$. This fact that is implicit in the work of Knaust,Odell, and Schlumprecht in \cite{KOS99} was mentioned in \cite{BKL10}, \cite{Lancien_course}, and the argument is detailed in \cite{Netillard22} (see also \cite{CauseyLancien23_CM}). Thus it immediately follows that $X$ is asymptotically B-$(2, \delta)$-convex where $2-\delta>2^{1/p}$.
\end{proof}

\begin{rema} In the local theory, if a Banach space $X$ is uniformly non-square (equivalently B-$(2, \delta)$-convex for some $\delta >0$), then $X$ is reflexive (see \cite{James64}). This is no longer true in the asymptotic setting since there are non-reflexive asymptotically uniformly smooth spaces, e.g. $c_0$.
\end{rema}

\begin{rema}It was shown in \cite{Giesy66} that a Banach space $X$ is B-convex if and only if its bidual $X^{**}$ is B-convex. In the asymptotic setting, this duality property does not hold anymore since $c_0^{**}=\ell_\infty$. For a separable counter-example, consider $X=Z_{c_0}^*$, the dual of Lindenstrauss space $Z_{c_0}$ as defined in \cite{Lindenstrauss71}. Then, $X$ is separable and admits an equivalent norm that is asymptotically uniformly smooth (see \cite[Theorem 2.1]{CauseyLancien19}), and hence $X$ asymptotically B-convex by Proposition \ref{prop:AUS->B-convex}, but $X^{**} = X \oplus \ell_1$ contains $\ell_1$ (see \cite{Lindenstrauss71}).
\end{rema}

\section{Asymptotic infratype}
\label{sec:infratype}

We introduce and study an asymptotic notion of infratype from the point of view of asymptotic structure.

\begin{defi}
\label{def:infratype}
Let $p \in [1, \infty ]$. We say that a Banach space $X$ has \emph{asymptotic infratype} $p$ if there exists $C>0$ such that for every directed set $D$, every $k \in \bN$, every weakly null tree $\dtreek \subset B_X$, and every scalar sequence $(a_1, \dots, a_k) \in \bR^k$, there exists $t \in D^k$ such that
\begin{equation}
\min_{(\vep_j)_{j=1}^k \in \{ \pm 1\}^k} \bnorm{ \sum_{j=1}^k \vep_j a_j x_{t_{\rest j}} }^p \leq C^p \sum_{j=1}^k \abs{a_j}^p.
\end{equation}
\end{defi}

Akin to the local setting, it follows from the submultiplicativity of the asymptotic $B$-convexity parameter $\bbar_X$ that asymptotic $B$-convexity implies non-trivial asymptotic infratype. The proof is essentially Pisier's proof adjusted to our setting.

\begin{prop}
\label{prop:asymp-B-convex->nontrivial-asymp-infratype}
Let $X$ be a Banach space.
If $X$ is asymptotically $B$-convex, then  $X$ has non-trivial asymptotic infratype.
\end{prop}

\begin{proof}
Since $X$ is asymptotically $B$-convex, there is $l_0\ge 2$ such that $l_0^{-1}\le \bbar_X(l_0)<1$. Thus one can find $p'\in [1,\infty)$ such that 
$\bbar_X(l_0)l_0^{\frac{1}{p'}}=1$. Let $p\in(1,\infty]$ such that $\frac{1}{p}+\frac{1}{p'}=1$ and pick $q\in(1,p)$. Then, fix $k \in \bN$, $D$ a directed set, $\dtreek \subset B_X$ a weakly null tree, and $\abar:=(a_1, \dots, a_k) \in \bR^k$. We now implement a classical blocking trick, adjusted to our asymptotic setting, to derive the asymptotic infratype $q$ inequality. 
Consider the partition $(0,\norm{a}_{\ell_q^k}]=\sqcup_{r\ge 0} A_r$, where $A_r = (\norm{\abar}_{\ell_q^k}l_0^{-(r+1)/q}, \norm{a}_{\ell_q^k}l_0^{-r/q}]$, for all $r \geq 0$.
Let $I_r = \left\{ 1 \leq j \leq k \colon \abs{a_j}\in A_r \right\}$ and observe that the set $B = \{ r \in \bN \colon I_r \neq \varnothing \}$ is finite and $\{1,\dots,k\}=\sqcup_{r\in B} I_r$. Moreover, since $\norm{a}_{\ell_q^r}^q\ge \sum_{j\in I_r} \abs{a_j}^q \ge \abs{I_r} \cdot\norm{a}_{\ell_q^k}^q l_0^{-(r+1)}$ it follows that $\abs{I_r}\le l_0^{r+1}$. Since $B$ is finite and $\abs{ a_j }\le l_0^{-r/q}\norm{\abar}_{\ell_q^k}$, for every $j\in I_r$, it follows from Lemma \ref{lem:pruning1} that we can find $t \in D^k$ such that for all $r \in B$,
\begin{equation}
\min_{(\vep_j)_{j=1}^k \in \{ \pm 1\}^k} \bnorm{ \sum_{j \in I_r} \vep_j  a_j x_{t_{\rest j}} } \leq |I_r| \bbar_X(|I_r|) l_0^{-r/q}\norm{\abar}_{\ell_q^k}.
\end{equation}
Using the submultiplicativity of $k\mapsto \bbar_X(k)$ and the fact that $k\mapsto k\bbar_X(k)$ is non-decreasing we deduce that
\begin{align*}
\min_{(\vep_j)_{j=1}^k \in \{ \pm 1\}^k} \bnorm{ \sum_{j=1}^k \vep_j a_j x_{t_{\rest j}} }& \leq  \sum_{r \in B} \min_{(\vep_j)_{j=1}^k \in \{ \pm 1\}^k} \bnorm{ \sum_{j \in I_r} \vep_j a_j x_{t_{\rest j}} } \le  \sum_{r \in B} |I_r| \bbar_X(|I_r|) l_0^{-r/q}\norm{\abar}_{\ell_q^k} \\
& \le \sum_{r \in B} l_0^{r+1} \bbar_X(l_0^{r+1})l_0^{-r/q}\norm{\abar}_{\ell_q^k} \le \sum_{r \in B} l_0^{r+1} \bbar_X(l_0)^{r+1} l_0^{-r/q}\norm{\abar}_{\ell_q^k}\\
& \le \sum_{r \in B} l_0^{r+1} l_0^{-{(r+1)}/p'} l_0^{-r/q}\norm{\abar}_{\ell_q^k} \le l_0^{1/p}\sum_{r \in B} (l_0^{1/p-1/q})^r\norm{\abar}_{\ell_q^k}.
\end{align*}
Since $\frac{1}{p}-\frac{1}{q}<0$ by the choice of $q$, we have that $\sum_{r \in B} (l_0^{1/p-1/q})^r \le \sum_{r =0}^\infty (l_0^{1/p-1/q})^r = \frac{l_0^{1/q-1/p}}{l_0^{1/q-1/p}-1}<\infty$, and hence we can take $C=\frac{l_0^{1/q}}{l_0^{1/q-1/p}-1}$.
\end{proof}

\begin{rema}
\label{rem:infra}
The proof shows that if $p' \in [1, \infty )$ and $p$ is its conjugate exponent, then $X$ has asymptotic infratype $q$ for all $q<p$ whenever $\bbar_X(l_0) \le  l_0^{-1/p'}$ for some $l_0 \in \bN$.
\end{rema}

The next corollary contains the equivalence between $(i)$ and $(iii)$ in Theorem \ref{thmB:asymp-infratype<->asymp-type} together with some additional quantitative information.

\begin{coro}
\label{thm:asymp-B-convex-asymp-infratype}
Let $X$ be a Banach space.
$X$ is asymptotically B-convex if and only if $X$ has non-trivial asymptotic infratype. \\
Moreover, if we denote by $\bar{\pi}_X$ the supremum of all $p \geq 1$ such that $X$ has asymptotic infratype $p$, then
\begin{equation}
\label{eq:asym-B-convex-infra}
\bar{\pi}_X= \lim_{k\to\infty} \frac{\log(k)}{\log(k \bbar_X(k))}.
\end{equation}
\end{coro}

\begin{proof}
The necessary condition is Proposition \ref{prop:asymp-B-convex->nontrivial-asymp-infratype}.
If $X$ has asymptotic infratype $p$, then there exists $C>0$ such that for all $k \geq 2$, $\bbar_X(k)\le C k^{1/p-1}$ and hence $\lim_{k\to\infty}\bbar_X(k)=0$ whenever $p>1$. Thus, by Corollary \ref{cor:asymp-B-convex-ell1} $X$ is asymptotically $B$-convex.  As for equality \eqref{eq:asym-B-convex-infra}, since for all $k \geq 2$, $k\bbar_X(k) \leq Ck^{1/p}$ we have
\begin{equation}
 \liminf\limits_{k\to\infty} \frac{\log(k)}{\log(k \bbar_X(k))} \geq \liminf\limits_{k\to\infty} \frac{\log(k)}{\log(C)+\frac{1}{p}\log(k)} = p,
 \end{equation}
and hence  $\liminf\limits_{k\to\infty} \frac{\log(k)}{\log(k \bbar_X(k))} \geq \bar{\pi}_X$. \\
For the other inequality, let $\vep>0$ and denote $q_\vep := \bar{\pi}_X +\vep$. Since $X$ does not have asymptotic infratype $\bar{\pi}_X +\vep$, it follows from Remark \ref{rem:infra} that necessarily $\bbar_X(k) > k^{-(1-1/q_{2\vep})}$ for all $k \geq 2$. Therefore
\begin{equation}
\limsup\limits_{k\to\infty} \frac{\log(k)}{\log(k \bbar_X(k))} \leq \limsup\limits_{k\to\infty} \frac{\log(k)}{\frac{1}{q_{2\vep}} \log(k)}  \le \bar{\pi}_X+2\vep.
\end{equation}
Thus $\limsup\limits_{k\to\infty} \frac{\log(k)}{\log(k \bbar_X(k))} \leq \bar{\pi}_X$.
\end{proof}

It remains to establish the connection between asymptotic infratype and asymptotic type. For this purpose, we need to introduce an additional parameter. 

For $k \in \bN$, we will denote by $\bar{\nu}_X(k)$ the smallest possible constant $\nu \geq 0$ such that for every directed set $D$, every weakly null tree $\dtreek \subset B_X$, and every scalar sequence $\abar:=(a_1, \dots, a_k) \in \bR^k$, there exists $t \in D^k$ so that
\begin{equation}
\label{eq:nu}
\bE \bnorm{ \sum_{j=1}^k \vep_j a_j x_{t_{\rest j}} } \leq \nu k \norm{\abar}
_{\ell_\infty^k}.
\end{equation}

The sequence $k\mapsto \bar{\nu}_X(k)$ behaves very much like $k\mapsto \bbar_X(k)$.

\begin{lemm} 
\label{lem:nu-submul}
Let $X$ be a Banach space. Then,
\begin{enumerate}[(i)]
\item $k\mapsto k\cdot\bar{\nu}_X(k)$ is non-decreasing,
\item[]and
\item for all $k,l \in \bN$, $\bar{\nu}_{X}(kl) \leq \bar{\nu}_{X}(k) \bar{\nu}_{X}(l)$.
\end{enumerate}
\end{lemm}

\begin{proof}
The proof of the monotonicity property is similar to the proof of $(i)$ in Lemma \ref{lem:b-submul} and is omitted.
Even though the argument is quite similar to $(ii)$ in Lemma \ref{lem:b-submul}, we provide some details for the submultiplicativity property. Let $k,l \in \bN$, $\nu_k>\bar{\nu}_{X}(k), \nu_l>\bar{\nu}_{X}(l)$, $D$ a directed set, $(x_t)_{t\in D^{kl}\setminus\{\emptyset\}} \subset B_X$ a weakly null tree, and $\abar:=(a_1, \cdots, a_{kl}) \in \bR^{kl}$ a scalar sequence. After pruning and relabelling if needed, we can assume that for all $t \in D^{kl}$ and all $1\le i \le k$,
\begin{equation}
\label{eq:nu-submul}
\bE_{\vep\in \{\pm 1\}^l} \bnorm{\sum_{j=1}^l \vep_{j} a_{(i-1)l+j} x_{t_{\rest (i-1)l+j}} }\le l\nu_l \max_{1\le j \le l} \abs{a_{(i-1)l+j}}.
\end{equation}
For every $\vep\in\{\pm 1\}^{kl}$ we will construct a weakly null tree of size $k$ goes as follows. This time we let $D_X$ be a weak convex symmetric neighborhood basis of $0$ in $X$ directed by reverse inclusion.  
Let $U \in D_X$ and $\vep\in\{\pm 1\}^{kl}$. We can find $V^{(U)} \in D_X$ such that $\underbrace{V^{(U)}+\dots+V^{(U)}}_{l~\mathrm{ times}} \subset U$. We can pick $t^{(U)}\in D^l$ such that $a_j x_{t^{(U)}_{\rest j}} \in V^{(U)}$ for all $1 \leq j \leq l$. Note that $t^{(U)}$ does not depend on $\vep$. Then, we let $y_U(\vep) := \sum_{j=1}^l \vep_j a_j x_{t^{(U)}_{\rest j}} \in U$, and observe that $(y_U(\vep))_{U \in D_X}$ is weakly null. Iterating we can find for all $1\le i \le k$ and $(U_1,U_2,\dots, U_i)\in D_X^i$, a node $t^{(U_1,U_2, \dots,U_i)}\in D^{il}$ such that 
\begin{enumerate}
    \item $a_j x_{ t^{(U_1,\dots,U_i)}_{\rest (i-1)l+j} } \in V^{(U_i)}$ for all $1 \leq j \leq l$, where $V^{(U_i)} \in D_X$ is such that $\underbrace{V^{(U)}+\dots+V^{(U)}}_{l~\mathrm{ times}} \subset U_i$
    \item $t^{(U_1,\dots,U_{i-1})}\prec t^{(U_1,\dots,U_i)}$.
\end{enumerate}
If for all $1\le i \le k$ we define 
\begin{equation}
\label{eq2:nu-submul}
y_{(U_1,\dots,U_i)}(\vep) := \sum_{j=1}^l \vep_{(i-1)l+j} a_{(i-1)l+j} x_{t^{(U_1,U_2, \dots,U_i)}_{\rest (i-1)l+j}},
\end{equation}
then for all $\vep\in \{-1,1\}^{kl}$, the tree $(\frac{ y_{u}(\vep) }{ \norm{ y_{u}(\vep) } })_{u \in D_X^{\le k}\setminus \{\emptyset\}}$ is weakly null and in $B_X$.

By definition of $\bar{\nu}_X(k)$, and after pruning one more time, there exists $u \in D_X^k$ such that for all $\vep\in \{-1,1\}^{kl}$, one has 
\begin{equation}
\label{eq3:nu-submul}
    \bE_{\eta} \bnorm{ \sum_{i=1}^k \eta_i y_{ u_{\rest i} }(\vep) } \leq \nu_k k \max_{1\le i \le k} \norm{y_{u_{\rest i}}(\vep)} .
\end{equation} 
Therefore, since for all $1\le i \le k$, we have 
\begin{equation*}
    \bE_{\vep\in\{\pm 1\}^{kl}}\norm{y_{u_{\rest i}}(\vep)} \stackrel{\eqref{eq2:nu-submul}}{=} \bE_{\vep\in \{\pm 1\}^l} \bnorm{\sum_{j=1}^l \vep_{j} a_{(i-1)l+j} x_{t^{(U_1,U_2, \dots,U_i)}_{\rest (i-1)l+j}} } \stackrel{\eqref{eq:nu-submul}}{\le} l \nu_l \ds\max_{1\le j \le l} \abs{a_{(i-1)l+j}},
\end{equation*}
it follows from \eqref{eq3:nu-submul} that 
\begin{align*}
    \bE_\vep \bE_{\eta} \bnorm{ \sum_{i=1}^k \eta_i y_{u_{\rest i} }(\vep) } & \le \nu_k k \max_{1\le i \le k} \bE_\vep\norm{y_{u_{\rest i}}(\vep)} \\
    & \le \nu_k k \max_{1\le i \le k}( l \nu_l \max_{1\le j \le l} \abs{a_{(i-1)l+j}}) = \nu_k\nu_l kl\max_{1\le i \le kl}\abs{a_i}.
\end{align*}
On the other hand, 
\begin{align*}
        \bE_\vep \bE_{\eta} \bnorm{ \sum_{i=1}^k \eta_i 
 y_{ u_{\rest i} }(\vep) } & = \bE_\vep \bE_{\eta} \bnorm{ \sum_{i=1}^k \eta_i 
 \sum_{j=1}^l \vep_{(i-1)l+j} a_{(i-1)l+j} x_{t^{(U_1,U_2, \dots,U_i)}_{\rest (i-1)l+j}} }\\
    & =  \bE_{\eta} \bE_\vep\bnorm{ \sum_{i=1}^k 
 \sum_{j=1}^l  \eta_i\vep_{(i-1)l+j} a_{(i-1)l+j} x_{t^{(U_1,U_2, \dots,U_i)}_{\rest (i-1)l+j}} } = \bE_\vep \bnorm{ \sum_{i=1}^k 
 \sum_{j=1}^l  \vep_{(i-1)l+j} a_{(i-1)l+j} x_{t^{(U_1,U_2, \dots,U_i)}_{\rest (i-1)l+j}} } \\
 & = \bE_\vep \bnorm{ \sum_{i=1}^{kl} \vep_{i} a_i x_{t^{(U_1,U_2, \dots,U_{k})}_{\rest i}} }.
\end{align*}
This means that there is $t\in D^{kl}$ such that 
\begin{equation*}
    \bE_\vep \bnorm{ \sum_{i=1}^{kl} \vep_{i} a_i x_{t_{\rest i}} } \le \nu_k\nu_l kl \norm{\abar}_{\ell_\infty^{kl}},
\end{equation*}
from which it follows that $\bar{\nu}_X(kl)\le \bar{\nu}_X(l)\bar{\nu}_X(k)$.
\end{proof}

In the following lemma, we compare the two parameters that we have introduced so far.
\begin{lemm}
\label{lem:b-mu-nu}
Let $X$ be a Banach space.
\begin{enumerate}[(i)]
\item For all $k \in \bN$, $0\le \bbar_X(k) \leq \bar{\nu}_{X}(k) \leq 1$.
\item For all $k \in \bN$, $\bbar_{X}(k)=1 \Longleftrightarrow \bar{\nu}_{X}(k)=1$.
\end{enumerate}
\end{lemm}

\begin{proof}
The first assertion simply follows by taking $\abar=(1,1,\dots,1)$ in \eqref{eq:nu}. 
For the second assertion, we just need to show that if $\bar{\nu}_{X}(k)=1$, then $\bbar_X(k)=1$. \\
We claim that for all $\nu<\bar{\nu}_k(X)$, we have $2^{k-1} (\nu-1) \leq \bbar_{X}(k)-1$, from which the conclusion immediately follows. To prove the claim, let $\nu < \bar{\nu}_{X}(k)$. It follows from the definitions that there exist a directed set $D$, a weakly null tree $\dtreek \subset B_X$, and a scalar sequence $0\neq \abar:=(a_1, \dots, a_k) \in \bR^k$ such that  for all $t\in D^k$,
\begin{equation}
\label{eq:13}
 \bE_\vep \bnorm{ \sum_{j=1}^k \vep_j a_j x_{t_{\rest j}} } > \nu k \norm{\abar}_{\ell_\infty^k}.
\end{equation}
and $\tau \in D^k$ such that $\min\limits_{(\vep_j)_{j=1}^k \in \{\pm 1 \}^k } \bnorm{\sum_{j=1}^k \vep_j x_{\tau_{\rest j}} } \leq \bbar_{X}(k) k$. We first obtain a crude upper bound, but sufficient for our purpose, between Rademacher averages and minima over sign choices. Given $y_1, \dots, y_k \in S_X$, let $\tilde{\vep}:=(\tilde{\vep}_j)_{j=1}^k \in \{-1, 1 \}^k$ such that
\begin{equation}
\min\limits_{(\vep_j)_{j=1}^k \in \{ \pm 1 \}^k } \Big\| \sum_{j=1}^k \vep_j y_j \Big\| = \Big\| \sum_{j=1}^k \tilde{\vep}_j y_j \Big\|.
\end{equation}
Then, 
\begin{align*}
\bE_\vep \bnorm{\sum_{j=1}^k \vep_j y_j} &= \frac{1}{2^k} \left( \sum_{(\vep_j)_{j=1}^k \in \{ \pm 1 \}^k \setminus \{ (\tilde{\vep}, -\tilde{\vep} \} } \bnorm{\sum_{j=1}^k \vep_j y_j } + 2 \min\limits_{(\vep_j)_{j=1}^k \in \{ \pm 1 \}^k } \bnorm{\sum_{j=1}^k \vep_j y_j } \right) \\
& \leq \frac{1}{2^k} \left((2^k-2)k + 2 \min\limits_{(\vep_j)_{j=1}^k \in \{ \pm 1 \}^k } \bnorm{ \sum_{j=1}^k \vep_j y_j }\right).
\end{align*}
By taking $y_j=x_{\tau_{\rest j}}$, for $1 \leq j \leq k$, in the inequality above, we have
\begin{equation}
\label{eq:15}
 \bE_\vep \bnorm{\sum_{j=1}^k \vep_j x_{\tau_{\rest j}} }\leq \left( \frac{(2^{k-1}-1)k}{2^{k-1}} + \frac{\min\limits_{(\vep_j)_{j=1}^k \in \{ \pm 1 \}^k } \bnorm{\sum_{j=1}^k \vep_j x_{\tau_{\rest j}}}}{2^{k-1}} \right) \leq k \left( \frac{2^{k-1}-1 + \bbar_{X}(k)}{2^{k-1}} \right).
\end{equation}
To get rid of the coefficients in \eqref{eq:13} and exploit \eqref{eq:15}, we invoke the contraction principle which tells us that 
\begin{equation}
\bE_\vep \bnorm{\sum_{j=1}^k \vep_j a_j x_{\tau_{\rest j}}}  \le \norm{\abar}_{\ell_\infty^k} \bE_\vep \bnorm{\sum_{j=1}^k \vep_j x_{\tau_{\rest j}}}.
\end{equation}
Therefore, 
\begin{equation}
\nu k \norm{\abar}_{\ell_\infty^k} \le k\norm{\abar}_{\ell_\infty^k} \left( \frac{\bbar_{X}(k) + 2^{k-1}-1}{2^{k-1}} \right).
\end{equation}
Since $\norm{\abar}_{\ell_\infty^k}\neq 0$, it follows that 
\begin{equation}
\nu  \le  \frac{\bbar_{X}(k) + 2^{k-1}-1}{2^{k-1}} ,
\end{equation}
which after some rearranging gives the claim.

\end{proof}

It is immediate that non-trivial asymptotic type implies non-trivial asymptotic infratype. It follows from Lemma \ref{lem:b-mu-nu} that the converse also holds.

\begin{theo} Let $X$ be a Banach space.
\label{thm:asymp-infratype->asymp-type}
If $X$ has non-trivial  asymptotic infratype, then $X$ has non-trivial asymptotic type.
\end{theo}

\begin{proof}
Assume that $X$ has non-trivial asymptotic infratype. Then, $X$ is asymptotic $B$-convex and thus $l_0^{-1}\le \bbar_X(l_0)<1$ for some $l_0\in \bN$. It follows from Lemma \ref{lem:b-mu-nu} that $l_0^{-1}\le \bar{\nu}_X(l_0)<1$. At this point, we can argue essentially as in the proof of Proposition \ref{prop:asymp-B-convex->nontrivial-asymp-infratype} to show that $X$ has non-trivial asymptotic type. Here are the details. Choose $p'\in [1,\infty)$ such that $\bar{\nu}_X(l_0)l_0^{\frac{1}{p'}} = 1$. Let $p\in(1,\infty]$ such that $\frac{1}{p}+\frac{1}{p'}=1$ and we pick $q\in(1,p)$. Now we fix $k \in \bN$, $D$ a directed set, $\dtreek \subset B_X$ a weakly null tree, and $\abar:=(a_1, \dots, a_k) \in \bR^k$ a scalar sequence. Consider the partition $(0,\norm{\abar}_{\ell_q^k}]=\sqcup_{r\ge 0} A_r$, where $A_r = (\norm{\abar}_{\ell_q^k}l_0^{-(r+1)/q}, \norm{\abar}_{\ell_q^k}l_0^{-r/q}]$, for all $r \geq 0$.
Let $I_r = \left\{ 1 \leq j \leq k \colon \abs{a_j}\in A_r \right\}$ and observe that the set $B = \{ r \in \bN \colon I_r \neq \varnothing \}$ is finite and $\{1,\dots,k\}=\sqcup_{r\in B} I_r$. Moreover, since $\norm{\abar}_{\ell_q^r}^q\ge \sum_{j\in I_r} \abs{a_j}^q \ge \abs{I_r} \cdot\norm{\abar}_{\ell_q^k}^q l_0^{-(r+1)}$ it follows that $\abs{I_r}\le l_0^{r+1}$. Since $B$ is finite, it follows from the pruning lemma that we can find $t \in D^k$ such that for all $r \in B$,
\begin{equation}
\bE_\vep\bnorm{\sum_{j\in I_r} \vep_j a_j x_{t_{\rest j}}} = \bE_\vep \bnorm{\sum_{j\in I_r} \vep_j a_j x_{t_{\rest j}}}  \leq |I_r| \bar{\nu}_X(|I_r|) \max_{j\in I_r}\abs{a_j} \le |I_r| \bar{\nu}_X(|I_r|)\norm{\abar}_{\ell_q^k}l_0^{-r/q}.
\end{equation}
Using the submultiplicativity of $k\mapsto \bar{\nu}_X(k)$ and the fact that $k\mapsto k\cdot \bar{\nu}_X(k)$ is non-decreasing we deduce that
\begin{align*}
\bE_\vep\bnorm{\sum_{j =1}^k \vep_j a_j x_{t_{\rest j}}} & \leq  \sum_{r \in B} \bE_\vep\bnorm{\sum_{j\in I_r} \vep_j a_j x_{t_{\rest j}}} \le  \sum_{r \in B} |I_r| \bar{\nu}_X(|I_r|)\norm{\abar}_{\ell_q^k}l_0^{-r/q}\\
& \le \sum_{r \in B} l_0^{r+1} \bar{\nu}_X(l_0^{r+1})l_0^{-r/q}\norm{\abar}_{\ell_q^k} \le \sum_{r \in B} l_0^{r+1} \bar{\nu}_X(l_0)^{r+1} l_0^{-r/q}\norm{\abar}_{\ell_q^k}\\
& \le \sum_{r \in B} l_0^{r+1} l_0^{-{(r+1)}/p'} l_0^{-r/q}\norm{\abar}_{\ell_q^k} \le l_0^{1/p}\sum_{r \in B} (l_0^{1/p-1/q})^r\norm{\abar}_{\ell_q^k}.
\end{align*}
Since $\frac{1}{p}-\frac{1}{q}<0$ by the choice of $q$, we have that $\sum_{r \in B} (l_0^{1/p-1/q})^r \le \sum_{r =0}^\infty (l_0^{1/p-1/q})^r = \frac{l_0^{1/q-1/p}}{l_0^{1/q-1/p}-1}<\infty$. It follows from the Kahane inequalities that $X$ has asymptotic type $q$.
\end{proof}

\begin{rema}
\label{rem:type}
The proof shows that if $p' \in [1, \infty )$ and $p$ is its conjugate exponent, then $X$ has asymptotic type $q$ for all $q<p$ whenever $\bar{\nu}_X(l_0) = l_0^{-1/p'}$ for some $l_0 \in \bN$.
\end{rema}

Combining Proposition \ref{prop:asymp-B-convex->nontrivial-asymp-infratype} and Theorem \ref{thm:asymp-infratype->asymp-type} we obtain the remaining equivalence in Theorem \ref{thmB:asymp-infratype<->asymp-type}.

\begin{coro}
Let $X$ be a Banach space. 
$X$ has non-trivial asymptotic type if and only if $X$ has non-trivial asymptotic infratype.
\end{coro}

\section{Asymptotic stable type}
\label{sec:stable-type}

The classical local notion of stable type involves $p$-stable random variables. We briefly recall their definition and some of their basic properties. If $\xi$ is a random variable defined on some probability space $(\Omega,\bP)$, its characteristic function is the map $\varphi_\xi\colon \bR \to \bC$ defined as $\varphi_\xi(t):=\bE(e^{it\xi})=\int_\Omega e^{it\omega}d\bP_\xi(\omega)$. Given $p\in(0,2]$, a random variable $\xi$ on a probability space is called \emph{$p$-stable} if its characteristic function is of the form $t \in \bR \mapsto \varphi_\xi(t)=e^{-c|t|^p}$, for some positive constant $c=c(p)$. Note that a centered Gaussian variable is $2$-stable. It is known that if $\xi$ is a $p$-stable random variable, then $\xi \notin L_p$ but $\xi \in L_q$ for all $q\in(0,p)$ (see \cite{AK06} for example). In fact, if $(\xi_i)_{i\in \bN}$ is a sequence i.i.d of $p$-stable random variables, then for all $q<p$ and $(a_i)_{i=1}^k\subset \bR$ 
\begin{equation}
\bnorm{\sum_{i=1}^k a_i \xi_i}_q= \norm{\xi_1}_q \Big(\sum_{i=1}^k \abs{a_i}^p\Big)^{1/p}.
\end{equation}

\begin{defi}
\label{def:asymp-stable-type}
Let $p \in (0,2]$. We say that a Banach space $X$ has \emph{asymptotic stable type $p$} if there exists $T_s>0$ such that for every sequence $(\xi_n)_{n\in\bN}$ of i.i.d. $p$-stable random variables, every $k \in \bN$, every directed set $D$, every weakly null tree $\dtreek \subset B_X$, and every scalar sequence $(a_1, \dots, a_k) \in \bR^k$, there is $t \in D^k$ such that
\begin{equation}
\bE \left( \bnorm{\sum_{j=1}^k \xi_j a_j x_{t_{\rest j}}}^{\frac{p}{2}} \right)^{2/p} \leq T_s \left(\sum_{j=1}^k \abs{a_j}^p\right)^{1/p}.
\end{equation}
\end{defi}

The proofs of the local analogs of the following statements can be implemented almost without any change in the asymptotic setting. We will thus only sketch them.

\begin{prop}
\label{prop:stable-type}
Let $p \in [1, 2]$ and $X$ be a Banach space.
\begin{enumerate}[(i)]
    \item If $X$ has asymptotic stable type $p$, then $X$ has asymptotic type $p$.
    \item If $X$ has asymptotic type $p$, then $X$ has asymptotic stable type $q$ for all $q \in (0,p)$.
    \item If $X$ has asymptotic stable type $1$, then $\ell_1$ is not asymptotically finitely representable in $X$.
\end{enumerate}
\end{prop}

\begin{proof}[Sketch of proof]
The first assertion is deducted from two well-known sets of inequalities. The first set are the Hoffmann-J\o rgensen inequalities \cite{HJ74}] which say that if $X$ is a Banach space and $(\xi_n)_{n\in\bN}$ a sequence of i.i.d. $s$-stable random variables with $s\in (0,2]$, then for each $p,q \in (0, s)$ if $s<2$, and each $p,q \in (0, \infty )$ if $s=2$, there exists a constant $K>0$ such that, for all $x_1,\dots, x_k$ in $X$,
\begin{equation}
\left( \bE \bnorm{\sum_{n=1}^k \xi_n x_n}^q \right)^{1/q} \leq K \left( \bE\bnorm{\sum_{n=1}^k \xi_n x_n }^p \right)^{1/p}.
\end{equation}
The other set of inequalities, which can be found in \cite{Pisier73}, states that if $(\varphi_n)_{n\in\bN}$ is a symmetric sequence of integrable real random variables and $(x_n)_{n=1}^k$ be a finite sequence in a Banach space $X$, then, for all $p \in [1, \infty )$,
\begin{equation*}
    \inf_{n\in \bN} \norm{\varphi_n}^p_{L_1} \bE_\vep\bnorm{ \sum_{n=1}^k \vep_n x_n }^p  \leq \bE \bnorm{ \sum_{n=1}^k \varphi_n x_n}^p. 
\end{equation*}

The second assertion follows essentially from Kahane inequalities and an inequality from \cite{Pisier73} which says that given $\abar:=(a_n)_{n \in \bN}$ a sequence of real numbers, $0 < r < q < p \leq 2$, and $(\xi_n)_{n\in \bN}$ a sequence of i.i.d. $q$-stable random variables, then for any $p$-stable random variable $g$, one has
\begin{equation}
    \norm{\xi_1}_{L_r} \norm{\abar}_{\ell_q} \leq \bnorm{ \Big(\sum_{n=1}^\infty \abs{a_n \xi_n}^p \Big)^{1/p}}_{L_r} \leq \frac{\norm{\xi_1}_{L_r} \norm{g}_{L_q}}{\norm{g}_{L_r}} \norm{\abar}_{\ell_q}.
\end{equation}

For the last assertion, assume that $\ell_1$ is asymptotically finitely representable in $X$ and that $X$ has asymptotic stable type $1$ with constant $C$.
Let $\zbar:=(z_n)_{n\in\bN} \in \ell_1$. In order to get a contradiction, let us prove that $\sum_{n} \abs{z_n \xi_n}$ converges a.s. \\
By our assumption, for every $k\ge 1$, there exists a weakly null tree $\dtreek \subset B_X$ such that 
\begin{equation}
\label{eq:ascv}
    \frac{1}{2} \norm{\abar}_{\ell_1^k} \leq \bnorm{\sum_{j=1}^k a_j x_{t_{\rest j}} } \leq \norm{\abar}_{\ell_1^k}
\end{equation}
for all $t\in D^k$ and all $\abar:=(a_1, \dots, a_k) \in \bR^k$. Since $X$ has asymptotic stable type $1$, there exists $t \in D^k$ such that
\begin{equation*}
    \left( \bE  \bnorm{ \sum_{j=1}^k \xi_j z_j x_{t_{\rest j}} }^{1/2} \right)^2 \leq C \norm{(z_1, \dots, z_k)}_{\ell_1^k},
\end{equation*}
and we get from the lower bound in \eqref{eq:ascv} that 
\begin{equation*}
    \left( \bE  \Big( \sum_{j=1}^k \abs{\xi_j z_j} \Big)^{\frac{1}{2}} \right)^2 \leq 2C \norm{\zbar}_{\ell_1}<\infty.
\end{equation*}
Therefore, $\sum_n \abs{z_n \xi_n}$ converges a.s. Since it is well known (see \cite{Pisier74} and references therein) that $\sum_n \abs{z_n \xi_n}$ converges a.s if and only if 
\begin{equation}
\label{eq2:ascv}
\sum_n \abs{z_n} \Big( 1 + \bigabs{\ln (\abs{z_n}^{-1})} \Big) < \infty,
\end{equation}
we get a contraction by considering $(z_n)_{n\in \bN} \in \ell_1$ that does not satisfy $\eqref{eq2:ascv}$.
\end{proof}

\section{Concluding remarks}

Certainly, more can be said about the asymptotic notions discussed in this note. We give one example. A Banach space $X$ has an unconditional asymptotic structure, if there exists a constant $C$ such that for every $k \in \bN$, every directed set $D$, every weakly null tree $\dtreek \subset B_X$, there is $t \in D^k$ so that for all $(\vep_1, \dots, \vep_k) \in \{-1, 1 \}^k$.
\begin{equation*}
    \bnorm{\sum_{j=1}^k \vep_j x_{t_{\rest j}} } \leq C \bnorm{\sum_{j=1}^k x_{t_{\rest j}} }.
\end{equation*} 

It follows that if a Banach space with an unconditional asymptotic structure has asymptotic Rademacher type $p \in (1, \infty)$, then $X$ has upper $\ell_p$ tree estimates, i.e., $X$ is in the class $\mathsf{A}_p$ in Causey's terminology (see \cite{Causey21} or \cite{CFL23}). In particular, $X$ is $q$-asymptotically uniformly smoothifiable for all $q \in (1,p)$. This is similar to the well-known fact that a Banach lattice (or with an unconditional basis) has type $p$, then it is $q$-convex for all $q \in (1,p)$ (see \cite[Corollary 1.f.9]{LT2}).

There are many natural questions that we did not touch upon in this article. The most natural ones are of a quantitative nature. The exact relationship between asymptotic infratype and asymptotic type must be elucidated. Quite remarkably, Talagrand has shown that when $p\in[1,2)$, type $p$ and infratype $p$ coincide \cite{Talagrand92} but the two notions differ when $p=2$ \cite{Talagrand04}. Whether these conclusions still hold in the asymptotic setting is certainly worth studying. While some questions could be swiftly answered, others might potentially be much more challenging.

In this note, we focused on the notion of type but a similar study should be carried out for the dual notion of cotype. The notion of asymptotic cotype $q$ for $q<\infty$ was also introduced in \cite{CDK19} where an asymptotic Maurey-Pisier $\ell_\infty$-theorem was proved. One can certainly define the notion of asymptotic supracotype $q$ where we reverse inequality \eqref{eq:infratype} and change the $\min$ to a $\max$. The local analog of this notion has also been studied via equivalent properties with different names such as the Orlicz property (when $q=2$) or Banach spaces whose identity operator is a $(q,1)$-summing map. Here again, Talagrand showed that cotype $q$ and supracotype $q$ coincide \cite{Talagrand92_supracotype}, but yet again the two notions differ when $q=2$ \cite{Talagrand94}. All these results should be examined from the point of view of asymptotic structure.

Developing asymptotic analogs of local notions from the point of view of asymptotic structure and understanding the similarities and discrepancies between the structural results about these linear notions is fundamental to further our understanding of the local theory as well as to uncover new applications that the local theory cannot settle. Some of these applications can be found in the nonlinear geometry of Banach space. It is worth noting that Bourgain, Milman, and Wolfson \cite{BMW86} showed that a Banach space $X$ has nontrivial (Rademacher) type if and only if it $\sup_{k\in\bN}\cdist{X}(\mathsf{H_k})=\infty$, where $\cdist{X}(\mathsf{H_k})$ denotes the smallest bi-Lipschitz distortion with which the $k$-dimensional Hamming cube can be embedded into $X$. Since $\co$ has non-trivial asymptotic type and $\cdist{\co}(M) \le 2$ for every separable metric space $M$ (this is due to Aharoni \cite{Aharoni74} with a suboptimal distortion and the optimal distortion $2$ was obtained by Kalton and Lancien in \cite{KL2008}), there is no chance to prove an asymptotic analog, for separable spaces, of the Bourgain-Milman-Wolfson characterization for the notion of asymptotic type studied in this note. Consequently, while asymptotic versions of local properties from the point of view of asymptotic structure might provide a nice, despite limited, analogy with the linear local theory, they might be ill-defined for nonlinear purposes, or that additional restrictions must be enforced. This means that other ways of defining asymptotic analogs of local properties, that are better suited for nonlinear problems, are needed. One possible approach is developed in \cite{BaudierFovelle} where a different asymptotic notion of type is introduced in order to understand the geometry of the countably branching Hamming cubes.

\begin{thank}
    This work was initiated while the second named
author was visiting Texas A$\&$M University in College Station. She wishes to thank Florent Baudier for his warm
hospitality and the excellent working environment and also Gilles Lancien for enlightening discussions.
\end{thank}

\bibliographystyle{alpha}

\begin{bibdiv}
\begin{biblist}


\bib{Aharoni74}{article}{
	author = {Aharoni, Israel},
	journal = {Israel J. Math.},
	pages = {284--291},
	title = {Every separable metric space is {L}ipschitz equivalent to a subset of {$c^{+}_{0}$}},
	volume = {19},
	year = {1974}
 }


\bib{AK06}{book}{
	author = {Albiac, Fernando},
 	author = {Kalton, Nigel J.},
	publisher = {Springer, New York},
	series = {Graduate Texts in Mathematics},
	title = {Topics in {B}anach space theory},
	volume = {233},
	year = {2006}
 }

\bib{BaudierFovelle}{article}{
	author = {Baudier, F.},
 	author = {Fovelle, A.},
	journal = {in preparation},
	number = {},
	pages = {},
	title = {On the asymptotic Enflo problem},
	volume = {},
	year = {}
 }
 
\bib{BKL10}{article}{
	author = {Baudier, F.},
 	author = {Kalton, N. J.},
  	author = {Lancien, G.},
	journal = {Studia Math.},
	number = {1},
	pages = {73--94},
	title = {A new metric invariant for {B}anach spaces},
	volume = {199},
	year = {2010}
 }


\bib{Beauzamy76-type}{incollection}{
	author = {Beauzamy, B.},
	booktitle = {S\'{e}minaire {M}aurey-{S}chwartz (1975--1976): {E}spaces {$L^p$}, applications radonifiantes et g\'{e}om\'{e}trie des espaces de {B}anach},
	pages = {Exp. Nos. 6-7, 29},
	publisher = {\'{E}cole Polytech., Palaiseau},
	title = {Op\'{e}rateurs de type {R}ademacher entre espaces de {B}anach},
	year = {1976}
 }


\bib{Beck62}{article}{
    AUTHOR = {Beck, Anatole},
     TITLE = {A convexity condition in {B}anach spaces and the strong law of
              large numbers},
   JOURNAL = {Proc. Amer. Math. Soc.},
    VOLUME = {13},
      YEAR = {1962},
     PAGES = {329--334},
}


\bib{BMW86}{article}{
	author = {Bourgain, J.},
	author = {Milman, V.},
	author = {Wolfson, H.},
	journal = {Trans. Amer. Math. Soc.},
	number = {1},
	pages = {295--317},
	title = {On type of metric spaces},
	volume = {294},
	year = {1986}
 }



\bib{Causey21}{article}{
	author = {Causey, Ryan M.},
	journal = {Studia Math.},
	number = {2},
	pages = {155--212},
	title = {Three and a half asymptotic properties},
	volume = {257},
	year = {2021}
 }


\bib{CDK19}{article}{
	author = {Causey, Ryan M.},
 	author = {Draga, Szymon},
  	author = {Kochanek, Tomasz},
	journal = {Trans. Amer. Math. Soc.},
	number = {11},
	pages = {8173--8215},
	title = {Operator ideals and three-space properties of asymptotic ideal seminorms},
	volume = {371},
	year = {2019},
	bdsk-url-1 = {https://doi.org/10.1090/tran/7759}
 }


\bib{CFL23}{article}{
	author = {Causey, R. M.},
	author = {Fovelle, A.},
        author = {Lancien, G.},
	journal = {Trans. Amer. Math. Soc.},
	number = {3},
	pages = {1895--1928},
	title = {Asymptotic smoothness in {B}anach spaces, three-space properties and applications},
	volume = {376},
	year = {2023}
 }


\bib{CauseyLancien23_CM}{article}{
	author = {Causey, R. M.},
        author = {Lancien, G.},
        journal = {Colloq. Math.},
	number = {2},
	pages = {281--324},
	title = {Asymptotic smoothness and universality in {B}anach spaces},
	volume = {172},
	year = {2023}
 }


\bib{CauseyLancien19}{article}{
	author = {Causey, R. M.},
        author = {Lancien, G.},
	journal = {Studia Math.},
	number = {2},
	pages = {109--127},
	title = {Prescribed {S}zlenk index of separable {B}anach spaces},
	volume = {248},
	year = {2019}
 }


\bib{Giesy66}{article}{
	author = {Giesy, Daniel P.},
	journal = {Trans. Amer. Math. Soc.},
	pages = {114--146},
	title = {On a convexity condition in normed linear spaces},
	volume = {125},
	year = {1966},
 }


\bib{HJ74}{article}{
	author = {Hoffmann-J\o rgensen, J\o rgen},
	journal = {Studia Math.},
	pages = {159--186},
	title = {Sums of independent {B}anach space valued random variables},
	volume = {52},
	year = {1974}
 }


\bib{James64}{article}{
	author = {James, Robert C.},
	journal = {Ann. of Math. (2)},
	pages = {542--550},
	title = {Uniformly non-square {B}anach spaces},
	volume = {80},
	year = {1964}
 }


\bib{James74}{article}{
	author = {James, Robert C.},
	journal = {Israel J. Math.},
	pages = {145--155},
	title = {A nonreflexive {B}anach space that is uniformly nonoctahedral},
	volume = {18},
	year = {1974}
 }

\bib{KL2008}{article}{
   author={Kalton, N. J.},
   author={Lancien, G.},
   title={Best constants for Lipschitz embeddings of metric spaces into
   $c_0$},
   journal={Fund. Math.},
   volume={199},
   date={2008},
   number={3},
   pages={249--272},
}

\bib{KOS99}{article}{
	author = {Knaust, H.},
 	author = {Odell, E.},
  	author = {Schlumprecht, Th.},
	journal = {Positivity},
	number = {2},
	pages = {173--199},
	title = {On asymptotic structure, the {S}zlenk index and {UKK} properties in {B}anach spaces},
	volume = {3},
	year = {1999},
 }


\bib{Lancien_course}{incollection}{
	author = {Lancien, Gilles},
	booktitle = {Topics in functional and harmonic analysis},
	pages = {77--101},
	publisher = {Theta, Bucharest},
	series = {Theta Ser. Adv. Math.},
	title = {A short course on nonlinear geometry of {B}anach spaces},
	volume = {14},
	year = {2013}
 }



\bib{Lindenstrauss71}{article}{
	author = {Lindenstrauss, Joram},
	journal = {Israel J. Math.},
	pages = {279--284},
	title = {On {J}ames's paper ``{S}eparable conjugate spaces''},
	volume = {9},
	year = {1971}
 }


\bib{LT2}{book}{
	author = {Lindenstrauss, Joram},
	author = {Tzafriri, Lior},
	publisher = {Springer-Verlag, Berlin-New York},
	series = {Ergebnisse der Mathematik und ihrer Grenzgebiete [Results in Mathematics and Related Areas]},
	title = {Classical {B}anach spaces. {II}},
	volume = {97},
	year = {1979}
 }



\bib{MMTJ95}{incollection}{
	author = {Maurey, B.},
	author = {Milman, V. D.},
	author = {Tomczak-Jaegermann, N.},
	booktitle = {Geometric aspects of functional analysis ({I}srael, 1992--1994)},
	pages = {149--175},
	publisher = {Birkh\"{a}user, Basel},
	series = {Oper. Theory Adv. Appl.},
	title = {Asymptotic infinite-dimensional theory of {B}anach spaces},
	volume = {77},
	year = {1995}
 }


\bib{MaureyPisier76}{article}{
	author = {Maurey, Bernard},
	author = {Pisier, Gilles},
	journal = {Studia Math.},
	number = {1},
	pages = {45--90},
	title = {S\'{e}ries de variables al\'{e}atoires vectorielles ind\'{e}pendantes et propri\'{e}t\'{e}s g\'{e}om\'{e}triques des espaces de {B}anach},
	volume = {58},
	year = {1976}
 }

\bib{Netillard22}{article}{
	author = {Netillard, Francois},
	journal = {North-West. Eur. J. Math.},
	pages = {103--110, i},
	title = {Banach spaces with the {B}lum-{H}anson property},
	volume = {8},
	year = {2022}}


\bib{Pisier74}{incollection}{
	author = {Pisier, G.},
	booktitle = {S\'{e}minaire {M}aurey-{S}chwartz 1973--1974: {E}spaces {$L^p$}, applications radonifiantes et g\'{e}om\'{e}trie des espaces de {B}anach},
	pages = {Exp. No. 7, 19 pp. (errata, p. E.1)},
	publisher = {\'{E}cole Polytech., Paris},
	title = {Sur les espaces qui ne contiennent pas de {$l\sp{1}\sb{n}$} uniform\'{e}ment},
	year = {1974}
 }


\bib{Pisier73}{article}{
	author = {Pisier, Gilles},
	journal = {C. R. Acad. Sci. Paris S\'{e}r. A-B},
	pages = {A1673--A1676},
	title = {Type des espaces norm\'{e}s},
	volume = {276},
	year = {1973}
 }


\bib{Talagrand92}{article}{
	author = {Talagrand, Michel},
	journal = {Invent. Math.},
	number = {1},
	pages = {41--59},
	title = {Type, infratype and the {E}lton-{P}ajor theorem},
	volume = {107},
	year = {1992}
 }


\bib{Talagrand92_supracotype}{article}{
	author = {Talagrand, Michel},
	journal = {Invent. Math.},
	number = {3},
	pages = {545--556},
	title = {Cotype and {$(q,1)$}-summing norm in a {B}anach space},
	volume = {110},
	year = {1992}
 }


\bib{Talagrand94}{article}{
	author = {Talagrand, Michel},
	journal = {Israel J. Math.},
	number = {1-3},
	pages = {181--192},
	title = {Orlicz property and cotype in symmetric sequence spaces},
	volume = {87},
	year = {1994}
 }


\bib{Talagrand04}{article}{
	author = {Talagrand, Michel},
	journal = {Israel J. Math.},
	pages = {157--180},
	title = {Type and infratype in symmetric sequence spaces},
	volume = {143},
	year = {2004}}


\end{biblist}
\end{bibdiv}

\end{document}